\documentclass[11pt]{article}
\usepackage{geometry}                
\usepackage{float}
\usepackage{graphicx}
\usepackage{amsmath,amssymb}

\usepackage{notes}

\title{A Note on Special Polynomials and Minimal Surfaces}
\author{Peter Connor}
\begin{document}
\maketitle

\noindent {\sc Abstract. } {\footnotesize}
When using Traizet's regeneration technique to construct minimal surfaces, the simplest nontrivial configurations are given as the roots of polynomials that satisfy a hypergeometric differential equation.  We exhibit examples of simple minimal surfaces exhibiting the same behavior. 

\noindent
{\footnotesize %2000 MSC numbers
2000 \textit{Mathematics Subject Classification}.
Primary 53A10; Secondary 49Q05, 53C42.
}

\noindent
{\footnotesize %key words and phrases
\textit{Key words and phrases}. 
Minimal surface, balance equations.
}

\section{Introduction}
Traizet's regeneration technique has been used to construct many new examples of embedded minimal surfaces in Euclidean space.  This technique uses the Weierstrass Representation for minimal surfaces and solves the Period Problem using the implicit function theorem in the neighborhood of a singular minimal surface.  This limit can be realized conformally as a noded Riemann surface, with the location of the nodes solving a set of balance equations.  The typical situation involves gluing tiny catenoidal necks between horizontal planes or flat annuli, with the location of the necks satisfying a set of balance equations.  In each of these cases, the simplest configurations are given as the roots of polynomials that satisfy a hypergeometric differential equation.  There still isn't an explanation for this phenomenon.   

In 2002, Traizet \cite{tr2} constructed one-parameter families of singly periodic minimal surfaces of arbitrary genus in the quotient that limit as a foliation of $\mathbb{R}^3$ by parallel planes with tiny catenoidal necks between each level.  One such example is Rieman's minimal surface.  The simplest cases place the catenoidal necks at the roots of Chebyshev polynomials.

%In 2002, Traizet \cite{tr4} constructed one-parameter families of finite total curvature minimal surfaces of arbitrary %genus that limit as an $N$-sheeted horizontal plane with tiny catenoidal necks between consecutive levels.  The %simplest examples are with three levels ($N=3$), in which case the only examples are the Costa Hoffman Meeks %family.  The catenoidal necks are located at the roots of $f(z)=A z^n+B$.  See figure \ref{figure:CHM}

%\begin{figure}[h]
	%\centerline{ 
		%\includegraphics[height=2in]{CHM1}
		%\hspace{.5in}
		%\includegraphics[height=2in]{CHM2}
		%	}
	%\caption{Two examples from the Costa Hoffman Meeks family, images courtesy of M. Weber}
	%\label{figure:CHM}
%\end{figure}

In 2005, Traizet and Weber \cite{tw1} constructed one-parameter families of screw motion invariant minimal surfaces of arbitrary genus that are asymptotic to the helicoid.  In this case, helicoids are glued together instead of planes and catenoids.  The simplest cases place the axes of the helicoids at the roots of Hermite polynomials.

In 2012, Connor and Weber  \cite{cw1} constructed three-parameter families of embedded, doubly periodic minimal surfaces of arbitrary genus in the quotient that limit as a foliation of $\mathbb{R}^3$ by parallel planes with tiny catenoidal necks between each level.  The simplest cases here place the catenoidal necks at the roots of Legendre polynomials.  See figure \ref{figure:(1,n)}.

\begin{figure}[h]
	\centerline{ 
		\includegraphics[height=1.8in]{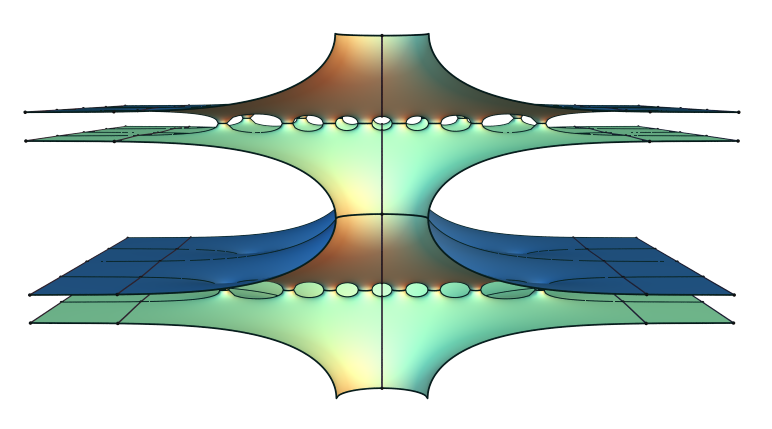}
		\hspace{.1in}
		\includegraphics[height=1.8in]{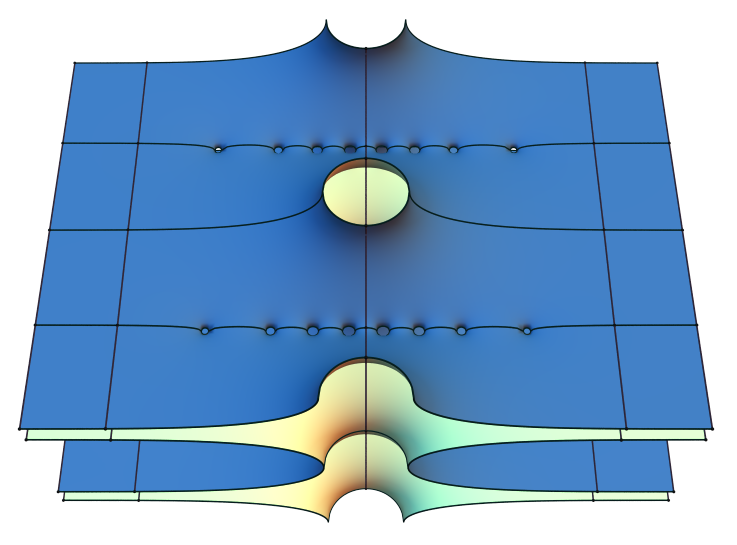}
			}
	\caption{Two views of a Connor-Weber surface}
	\label{figure:(1,n)}
\end{figure}

In 2012, Li  \cite{kli1} constructed one-parameter families of embedded, singly periodic surfaces with 6 Scherk ends of arbitrary genus that limit as three parallel planes connected by tiny catenoidal necks.  In this setting, the catenoidal necks are also located at the roots of polynomial solutions to a hypergeometric differential equation.

In this paper, we prove the existence of complete minimal surfaces with two horizontal annular ends and any number of catenoidal ends with vertical limiting normals such that the location of the catenoidal ends satisfy a similar set of balance equations as arise when one constructs surfaces using Traizet's regeneration technique.  This is done by considering the possible Weierstrass data for minimal surfaces with these properties, formulating the equations for the period problem, and then solving the period problem.  

\section{Weierstrass Representation}
We use the Weierstrass Representation of minimal surfaces, which may be written as
\[
X(z)=\re\int_{z_0}^z\left(\frac{1}{2}\left(\frac{1}{G}-G\right)dh,\frac{i}{2}\left(\frac{1}{G}+G\right)dh,dh\right)
\]
where $z,z_0\in\Sigma$ with $\Sigma$ a Riemann surface, $G$ is a meromorphic function called the Gauss map, and $dh$ is a meromorphic one-form called the height differential.  One issue is that $X$ depends on the path of integration.  

%The examples we are constructing are translation invariant in the direction of a vector $\Lambda$.  The map $X:\Sigma\rightarrow\R^3$ is well defined provided that 
%\[
%\re\int_\gamma\left(\frac{1}{2}\left(\frac{1}{G}-G\right)dh,\frac{i}{2}\left(\frac{1}{G}+G\right)dh,dh\right)=\vec{0}\mod(\Lambda)
%\]
%for all closed curves $\gamma$ in $\Sigma$.  This is called the period problem, and it can be expressed as

%\[
%\overline{\int_\gamma \frac{1}{G}dh}-\int_\gamma Gdh=0 \mod(\Lambda)
%\]
%and
%\[
%\re\int_\gamma dh=0 \mod(\Lambda)
%\]

%A second issue is that we want $X$ to be regular.  This is ensured by requiring that $G$ has either a zero or pole at $p\in\Sigma$ if and only if $dh$ has a zero at $p$ with the same multiplicity.

We want to construct a complete, singly periodic minimal surface with period $(0,\pi,0)$, two horizontal annular ends, and a sequence of catenoidal ends with vertical limiting normals.  Place the annular ends at $0$ and $\infty$, the catenoidal ends at $p_1, p_2,\ldots, p_n$, and assume that the vertical limiting normal vector points up at the catenoidal ends and down at the annular ends.  Then the domain of the map $X$ is $\Sigma=\C\setminus\{0,p_1,p_2,\ldots,p_n\}$.  Here are the conditions such that $(\Sigma, G, dh)$ provide the Weierstrass data for such a surface:
\begin{enumerate}
\item
The zeros of $dh$ are the zeros and poles of $G$, with the same multiplicity.  This ensures the map $X$ is regular.
\item
$G$ and $dh$ have zeros of order $k\geq 1$ and $k-1$ at $0$, respectively, and $G$ and $dh$ have zeros of order $j\geq 1$ and $j-1$, respectively, at $\infty$.  This ensures the ends at $0$ and $\infty$ are horizontal annular ends with downward pointing vertical limiting normal vectors.
\item
$G$ and $dh$ have simple poles at $p_j$ with residue $\alpha_j\in\R\setminus \{0\}$, $j=1,2,\ldots,n$.  This ensures that the ends at each $p_j$ are catenoidal ends with upward pointing vertical limiting normal vectors.  The residue $\alpha_j$ gives the neck size of the corresponding catenoidal end.
\item
\[
\re\int_{C(p_j)}\left(\frac{1}{2}\left(\frac{1}{G}-G\right)dh,\frac{i}{2}\left(\frac{1}{G}+G\right)dh,dh\right)=(0,0,0)
\]
where $C(p_j)$ is a simple closed curve around the puncture $p_j$, for $j=1,2,\ldots, n$.
\item
\[
\re\int_{C(0)}\left(\frac{1}{2}\left(\frac{1}{G}-G\right)dh,\frac{i}{2}\left(\frac{1}{G}+G\right)dh,dh\right)=(0,0,0)\;\text{mod}(0,\pi,0)
\]
and 
\[
\re\int_{C(\infty)}\left(\frac{1}{2}\left(\frac{1}{G}-G\right)dh,\frac{i}{2}\left(\frac{1}{G}+G\right)dh,dh\right)=(0,0,0)\;\text{mod}(0,\pi,0)
\]
where $C(0)$ and $C(\infty)$ are simple closed curves around $0$ and $\infty$.
\end{enumerate}
Together, $4$ and $5$ are referred to as the Period Problem.

Let
\[
G(z)=z\sum_{j=1}^n\frac{\alpha_j}{z-p_j}
\]
and 
\[
dh=\sum_{j=1}^n\frac{\alpha_j}{z-p_j}dz
\]
with 
\[
\sum_{j=1}^n\alpha_j=0
\]
Then all of the conditions on the Weierstrass data are automatically satisfied, except for the periods at each puncture $p_j$.  As we are computing integrals over domains in the complex plane, the period integrals can be calculated as residues.  The selection of the Weierstrass data reduces the period problem to 

\[
\res_{p_k}(Gdh)=0
\]
for $k=1,2,\ldots,n$, and 

\[
\res_{p_k}(Gdh)=\res_{p_k}\left(z\left(\sum_{j=1}^n\frac{\alpha_j}{z-p_j}\right)^2\right)=p_k \alpha_k \left(2\left(\sum_{j\neq k}\frac{\alpha_j}{p_k-p_j}\right)+ \frac{\alpha_k}{p_k}\right)
\]

Let 
\[
F_k=2\sum_{j\neq k}\frac{\alpha_j}{p_k-p_j}+\frac{\alpha_k}{p_k}
\]
If $F_k=0$ for $k=1,2,\ldots n$ then the period problem is solved.  The equations $F_k=0$ are referred to as the balance equations.  The Residue Theorem implies that one of the balance equations doesn't need to be solved.  The balance equations are equivalent to the balance equations in \cite{cw1} that produce doubly periodic examples with two levels in the quotient.  This isn't surprising because the Gauss map and height differential are defined similarly in both cases.

\section{Solving the balance equations}
Let's consider examples with one catenoid end opening downward and $n$ catenoid ends with the same necksize opening upward.  Then the balance equations are solved when the locations of the catenoidal ends are given by the roots of the polynomial
\[
f_n(z)=\sum_{j=0}^n {n\choose j}^2 z^j
\]
which is related to the Legendre polynomial $L_n(z)$ by the transformation
\[
f_n(z)=(1-z)^n L_n\left(\frac{1+z}{1-z}\right)
\] 
\begin{theorem}
Let $n\in\N$, $\alpha_j=1/n$ for $j=1,2,\ldots, n$, $\alpha_{n+1}=-1$, $p_{n+1}=1$, and $p_1,p_2,\ldots, p_n$ be the roots of the polynomial $f_n(z)$.  Then $F_k=0$ for $k=1,2,\ldots, n+1$, and $X(\Sigma)$ is a complete, immersed, singly periodic minimal surface with horizontal annular ends at $0$ and $\infty$, upward pointing catenoidal ends with vertical limiting normals at $p_1,p_2,\ldots,p_n$, and a horizontal downward pointing catenoidal end with vertical limiting normal at $p_{n+1}=1$.
\end{theorem} 

\begin{figure}[h]
	\centerline{ 
		\includegraphics[height=1.6in]{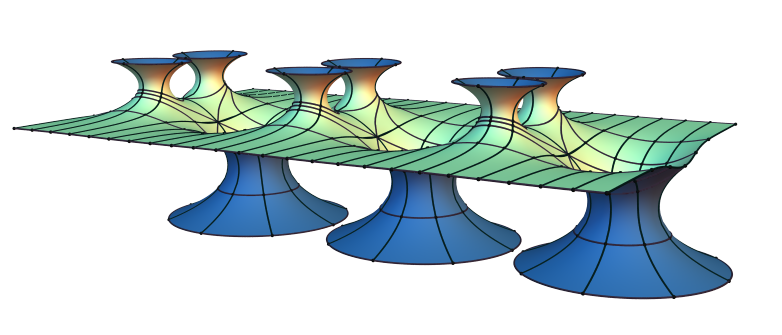}
		\hspace{.2in}
		\includegraphics[height=1.6in]{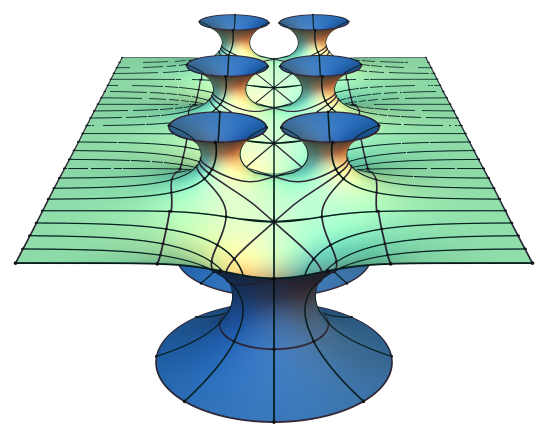}
			}
	\caption{Two views of a singly periodic minimal surface with three catenoidal ends in the quotient}
	\label{figure:(2,5)}
\end{figure}

\begin{proof}
In this case, the balance equations are
\[
F_k=\frac{2}{n}\sum_{j=1,j\neq k}^n\frac{1}{p_k-p_j}-\frac{2}{p_k-1}+\frac{1}{np_k}, \hspace{.1in} k=1,2,\ldots, n
\]
and 
\[
F_{n+1}=\frac{2}{n}\sum_{j=1}^n\frac{1}{p_k-p_j}-1
\]
Recall that it is sufficient to solve $F_k=0$ for $k=1,2,\ldots, n$.  

The polynomial $f_n$ satisfies the hypergeometric differential equation
\[
z(1-z)f_n''(z)+(1+(2n-1)z)f_n'(z)-n^2f_n(z)=0
\]
and thus has simple zeros.  Since $f_n$ has simple zeros, for each root $p_k$,
\[
f_n''(p_k)=2f_n'(p_k)\sum_{j=1,j\neq k}^n\frac{1}{p_k-p_j}
\]
Plugging into the above differential equation, we get that
\[
2p_k(1-p_k)\sum_{j=1,j\neq k}^n\frac{1}{p_k-p_j}+1+(2n-1)p_k=0
\]
which implies that $F_k$=0.
\end{proof}

\begin{figure}[h]
	\centerline{ 
		\includegraphics[height=1.6in]{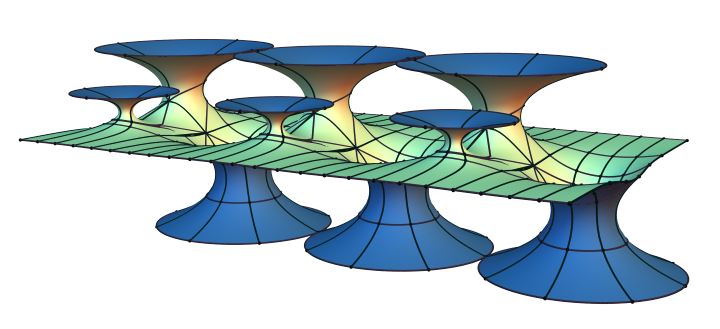}
		\hspace{.2in}
		\includegraphics[height=1.6in]{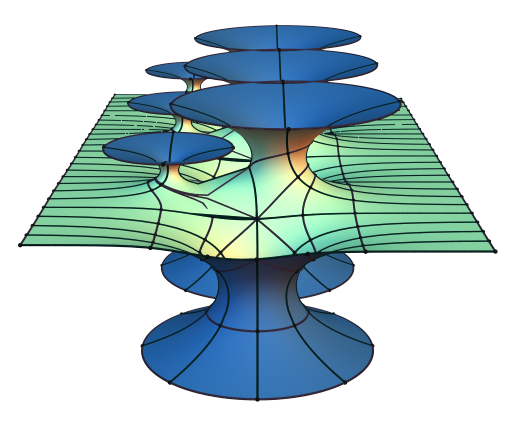}
			}
	\caption{Two views of a less symmetric example with three catenoidal ends in the quotient}
	\label{figure:nonsym}
\end{figure}

By adjusting the necksizes of the catenoids, one can construct less symmetric examples, even when there are only three catenoidal ends in the quotient.  The example of this type in figure \ref{figure:nonsym} has Weierstrass data
\[
G(z)=z\left(\frac{1}{4\left(z+\frac{11+3\sqrt{13}}{2}\right)}+\frac{3}{4\left(z-\frac{-7+\sqrt{13}}{6}\right)}-\frac{1}{z-1}\right)
\]
and 
\[
dh=\left(\frac{1}{4\left(z+\frac{11+3\sqrt{13}}{2}\right)}+\frac{3}{4\left(z-\frac{-7+\sqrt{13}}{6}\right)}-\frac{1}{z-1}\right)dz
\]

\addcontentsline{toc}{section}{References}
\bibliographystyle{plain}
\bibliography{minlit}

\label{sec:liter}

\end{document}